\documentclass[twoside,10pt]{article}
\usepackage{multicol}
\usepackage{graphicx}
\usepackage{xcolor}
\usepackage{amssymb}
\usepackage{amsmath,amsfonts}
\usepackage{amsthm}
\usepackage{enumerate}
\usepackage{psfrag}
\usepackage{array}
\input xy
\xyoption{all}
\usepackage[pagewise,displaymath, mathlines]{lineno}
\theoremstyle{plain}
\newtheorem{teo}{Theorem}[section]
\newtheorem{cor}[teo]{Corollary}

\newtheorem{prop}[teo]{Proposition}
\theoremstyle{definition}
\newtheorem{defi}[teo]{Definition}
\newtheorem{obs}[teo]{Remark}

\numberwithin{equation}{section}

\makeatletter
\def\blfootnote{\gdef\@thefnmark{}\@footnotetext}
\makeatother

\begin{document}
\thispagestyle{plain}

\bigskip

\begin{center}
{\huge{A note on homological systems}}
\end{center}

\begin{center}
{Jes\'us Efr\'en P\'erez Terrazas}
\end{center}

\begin{abstract}
We give an elementary short proof of the known fact that the
category $\mathfrak{F} \left( \Delta \right)$ of $\Delta -$filtered
modules, associated to a given finite homological system $\left(
\Delta ; \Omega, \leq \right),$ is closed under direct summands.
\end{abstract}

Homological systems are generalizations of the standard modules,
associated to quasi-hereditary algebras, and so they have similar
features (see \cite{R} and \cite{MSX}).

Homological systems in module categories over pre-ordered sets were
introduced in \cite{MSX}. There they show, see Corollary 3.16, that
the associated category of the $\Delta -$filtered modules is closed
under direct summands: if $M$ is a $\Delta -$filtered module and
there is a decomposition $M = M_1 \oplus M_2,$ then $M_1$ and $M_2$
are $\Delta -$filtered modules.

In this note we use a \emph{height map associated to a pre-ordered
set} defined below, a function determined by the relative position
of the elements of the pre-ordered set, to provide a direct proof of
that property.

\smallskip

From now on, $\Lambda $ is an Artin $R-$algebra, where $R$ is an
artinian commutative ring, $\Lambda -$mod is the full subcategory of
the finitely generated left $\Lambda-$modules, and all modules and
homomorphisms considered in this note belong to $\Lambda -$mod.

\begin{defi} \begin{upshape} A homological system $\left( \Delta ; \Omega , \leq
\right)$ consists of the following:

HS1: A finite pre-ordered set $\left( \Omega , \leq \right):$
$\Omega$ is a finite not empty set and $\leq $ is a reflexive and
transitive relation on $\Omega .$

HS2: The set $\Delta = \left\{ \Delta _{\omega} \right\}_{\omega \in
\Omega},$ where $\Delta _{\omega}$ is indecomposable in $\Lambda
-$mod and $\omega \neq \omega'$ implies $\Delta _{\omega} \ncong
\Delta _{\omega '}.$

HS3: If ${\rm Hom}_{\Lambda} \left( \Delta _{\omega} , \Delta
_{\omega '} \right) \neq 0$ then $\omega \leq \omega ' .$

HS4: If ${\rm Ext}_{\Lambda} \left( \Delta _{\omega} , \Delta
_{\omega '} \right) \neq 0$ then $\omega \leq \omega ' $ and it is
not true that $\omega ' \leq \omega .$
\end{upshape} \end{defi}

\begin{defi} \begin{upshape} Given the set $\Delta = \left\{ \Delta _{\omega} \right\}_{\omega \in
\Omega},$ as in HS2, denote by $\mathfrak{F} \left( \Delta \right)$
the full subcategory of $\Lambda -$mod of those $M$ having a $\Delta
-$filtration, that is a filtration of the form $$\{ 0 \} = M_0
\subsetneq M_1 \subsetneq \cdots \subsetneq M_t = M$$ such that
$M_{i} / M_{i-1} \cong \Delta _{\omega},$ for some $\omega$ that
depends of $i,$ and for each $i \in \left\{ 1, 2, \ldots , t
\right\}.$ So $\mathfrak{F} \left( \Delta \right)$ is closed under
isomorphisms, extensions and it contains zero objects. The modules
in $\mathfrak{F} \left( \Delta \right)$ are called $\Delta
-$filtered modules.
\end{upshape} \end{defi}

Given a pre-ordered set $\left( \Omega , \leq \right),$ the relation
$\omega \sim \omega'$ if and only if $\omega \leq \omega '$ and
$\omega ' \leq \omega$ is an equivalence relation, and so there is a
canonical surjective order-preserving function $\pi : \Omega
\longrightarrow \Omega / \sim ,$ where $\left( \Omega / \sim ,
\preceq \right)$ is a poset and $\preceq$ is the order induced by
$\leq .$

\begin{defi} \begin{upshape} Let $\left( \Omega ' , \preceq \right)$
be a finite poset. The height function $h' : \Omega '
\longrightarrow \mathbb{N}$ of the poset is defined recursively: the
elements of height 1 are the minimal elements of the poset, and the
elements of height $n+1$ are the minimal elements of the induced
poset on $\Omega - \left\{ \omega \in \Omega \; \mid \; h' \left(
\omega \right) \in \left\{ 1,2, \ldots , n \right\} \right\}.$ Given
a finite pre-ordered set $\left( \Omega, \leq \right),$ the
composition $h: \Omega \stackrel{\pi}{\longrightarrow} \Omega / \sim
\stackrel{h'}{\longrightarrow} \mathbb{N}$ is the height function of
the pre-ordered set.
\end{upshape} \end{defi}

\begin{obs} \begin{upshape} \label{immediate} Given a finite pre-ordered
set $\left( \Omega, \leq \right)$ and its height function $h,$ if
$\omega \leq \omega ' $ and $\omega ' \nleq \omega $ then $h \left(
\omega  \right) < h \left( \omega ' \right),$ and if $h \left(
\omega  \right) < h \left( \omega ' \right)$ then $\omega ' \nleq
\omega .$

So we get, for $\left( \Delta ; \Omega , \leq \right)$ a homological
system, that if $h \left( \omega \right)
> h \left( \omega ' \right)$ then ${\rm Hom}_{\Lambda} \left( \Delta
_{\omega} , \Delta _{\omega '} \right) = 0,$ and if $h \left( \omega
\right) \geq h \left( \omega ' \right)$ then ${\rm Ext}_{\Lambda}
\left( \Delta _{\omega} , \Delta _{\omega '} \right) = 0.$
\end{upshape} \end{obs}

\begin{defi} \begin{upshape} Let $F = \left\{ M_0 \subsetneq M_1 \subsetneq \cdots
\subsetneq M_t \right\}$ be a $\Delta -$filtration of $M \in
\mathfrak{F} \left( \Delta \right).$ We will denote by $\ell \left(
F \right) = t$ the length of $F,$ and by $\ell _{\omega} \left( F
\right)$ the number of factors isomorphic to $\Delta _{\omega},$ and
so $\ell \left( F \right) = \sum _{\omega \in \Omega} \ell _{\omega}
\left( F \right).$
\end{upshape} \end{defi}

\begin{prop} \label{arrangement} Let $\left( \Delta ; \Omega , \leq
\right)$ be a homological system, $M \in \mathfrak{F} \left( \Delta
\right) - \{ 0 \} ,$ and $F = \left\{ M_0 \subsetneq M_1 \subsetneq
\cdots \subsetneq M_t \right\}$ be a $\Delta -$filtration. Then
there exists a $\Delta -$filtration $F' = \left\{ M'_0 \subsetneq
M'_1 \subsetneq \cdots \subsetneq M'_t \right\}$ of $M$ with $\ell
_{\omega} \left( F \right) = \ell _{\omega} \left( F' \right),$ for
each $\omega \in \Omega ,$ and such that: $1 \leq i \leq j \leq t,$
$M'_i / M'_{i-1} \cong \Delta _{\omega}$ and $M'_j / M'_{j-1} \cong
\Delta _{\omega '},$ imply $h \left( \omega \right) \geq h \left(
\omega ' \right).$
\end{prop}

{\bf Proof:} With the notation of the statement, but starting with
the filtration $F,$ let us assume that for a fixed $i$ and $j = i+1$
we have $h \left( \omega \right) < h \left( \omega ' \right).$ Then,
from HS4, we have that the exact sequence $$\xymatrix{0 \ar[r] &
M_{i}/M_{i-1} \ar[r] & M_{i+1}/M_{i-1} \ar[r] & M_{i+1} / M_i \ar[r]
& 0}$$ splits, so there exists $\overline{N}$ submodule of
$M_{i+1}/M_{i-1}$ such that $\overline{N} \cong \Delta _{\omega '}$
and $\left( M_{i+1}/M_{i-1} \right) / \overline{N} \cong \Delta
_{\omega}.$

From the third theorem of isomorphism we get for $N,$ the pre-image
of $\overline{N}$ under the canonical epimorphism $p : M_{i+1}
\longrightarrow M_{i+1}/M_{i-1},$ that $F_1 = \left\{ M_0 \subsetneq
M_1 \subsetneq \cdots \subsetneq M_{i-1} \subsetneq N \subsetneq
M_{i+1} \subsetneq \cdots \subsetneq M_t \right\}$ is a $\Delta
 -$filtration with $\ell _{\omega} \left( F \right) = \ell _{\omega} \left( F_1 \right)$
 for each $\omega,$ satisfying the order condition in positions $i$ and $i+1.$

 We can repeat this process generating new $\Delta -$filtrations, and in
 a finite number of steps to obtain one as in the statement.
 $\hfill \square$

For the proof of the following proposition, let us recall that given
$L, N \in \Lambda -$Mod, the \emph{trace} ${\rm tr}_{L} \left( N
\right)$ of $L$ on $N$ is the sum of all the images of the
homomorphisms from $L$ to $N.$

\begin{prop} \label{slices} Let $\left( \Delta ; \Omega , \leq
\right)$ be a homological system and $M \in \mathfrak{F} \left(
\Delta \right).$ Let $\left\{ 1, 2, \ldots , a \right\}$ be the
image of the height function $h.$ There exists a filtration $\{ 0 \}
= W_{a+1} \subseteq W_{a} \subseteq \cdots \subseteq W_1 = M$ such
that
$$W_i / W_{i+1} \cong \displaystyle \bigoplus _{\substack{\omega \in \Omega \\ h \left( \omega \right) = i}} n_{\omega} \Delta _{\omega}$$
with $n_{\omega} \in \{ 0 \} \cup \mathbb{N}.$ Let us call a such
filtration an $h-$filtration. Also, the numbers $n'_{\omega}$ of any
other $h-$filtration are the same.
\end{prop}

{\bf Proof:} The existence of at least one $h-$filtration follows
directly from the Proposition \ref{arrangement} and HS4.

Now assume that $\{ 0 \} = W'_{a+1} \subseteq W'_{a} \subseteq
\cdots \subseteq W'_1 = M$ is another $h-$filtration. By HS3 and the
definition of trace we have
$$W_{a} = {\rm tr}_{W_{a}} \left( W_a \right) = {\rm tr}_{W_{a}} \left( M \right) = {\rm
tr}_{W_{a}} \left( W'_a \right) = W'_{a}$$ and so, by the
Krull-Schmidt-Remak Theorem, it follows that $n_{\omega} =
n'_{\omega}$ for any $\omega \in \Omega$ of height $a.$

We repeat the argument in the quotients $M/W_{a},$ $M/W_{a-1},$ ...,
$M /W_2$ in order to verify the statement. $\hfill \square$

\begin{cor} \label{additivity} Let $\left( \Delta ; \Omega , \leq
\right)$ be a homological system.
\begin{enumerate}
\item Let $M \in \mathfrak{F} \left( \Delta \right)$ and $F$ and
$F'$ be $\Delta -$filtrations of $M.$ Then $\ell _{\omega} \left( F
\right) = \ell _{\omega} \left( F' \right)$ for each $\omega \in
\Omega .$ It follows that the $\Delta-$length and the number of
$\Delta -$factors are well defined for $M.$
\item Let $L,N \in \mathfrak{F} \left( \Delta \right)$ and $\xymatrix{0 \ar[r] & L \ar[r] & M \ar[r] & N \ar[r] &
0}$ be an exact sequence. Then $\ell _{\omega} \left( M \right) =
\ell _{\omega} \left( L \right) + \ell _{\omega} \left( N \right)$
for each $\omega \in \Omega .$
\item $\mathfrak{F} \left( \Delta \right)$ is closed under direct
summands.
\end{enumerate}
\end{cor}

{\bf Proof:}

1.- If $F$ is a $\Delta -$filtration of $M,$ and $\underline{F}$ is
an $h-$filtration obtained from $F,$ as in the proofs of the
propositions \ref{arrangement} and \ref{slices}, then $\ell_{\omega}
\left( F \right) = n_{\omega}$ for each $\omega ,$ where
$n_{\omega}$ are the coefficients associated to $\underline{F}$ as
in \ref{slices}. Then the claim follows by the uniqueness of those
coefficients.

2.- It is a direct consequence of the previous item.

3.- Consider an $h-$filtration of $M$ and assume that $M = M_1
\oplus M_2.$

Then, by additivity of the trace, we have $W_a = {\rm tr}_{W_a}
\left( M \right) = {\rm tr}_{W_a} \left( M_1 \right) \oplus {\rm
tr}_{W_a} \left( M_2 \right),$ and so ${\rm tr}_{W_a} \left( M_j
\right)$ is a submodule of $M_j$ which is in $\mathfrak{F} \left(
\Delta \right),$ for $j \in \left\{ 1, 2 \right\}.$

We can repeat this argument for the quotients $M /W_a \cong \left(
M_1 / {\rm tr}_{W_a} \left( M_1 \right) \right) \oplus \left( M_2 /
{\rm tr}_{W_a} \left( M_2 \right) \right)$ in order to get bigger
submodules of $M_1$ and $M_2$ that are in $\mathfrak{F} \left(
\Delta \right).$

Repeating the procedure we obtain $h-$filtrations of $M_1$ and
$M_2,$ so they belong to $\mathfrak{F} \left( \Delta \right).$
$\hfill \square$

\bigskip

E. P\'erez

Facultad de Matem\'aticas

Universidad Aut\'onoma de Yucat\'an

M\'erida, M\'exico

jperezt@correo.uady.mx

\end{document}